\documentclass[12pt]{article}
\usepackage{natbib}
\usepackage{amssymb,amsmath,amsfonts,amssymb}
\usepackage{graphicx}
\usepackage{subcaption}
\usepackage{multicol}
\input{epsf.sty}
 
\usepackage{comment}

\input{epsf.sty}
\newcommand{\NI}{\noindent}
\newtheorem{theorem}{\NI{\bf Theorem}}[section]
\newtheorem{lemma}{\NI\bf Lemma}[section]
\newtheorem{prop}{\NI\bf Proposition}[section]
\newtheorem{cor}{\NI\bf Corollary}[section]

\newtheorem{defn}{\NI\bf Definition}[section]
\newtheorem{example}{\NI\bf Example}[section]
\newcommand{\bt}{\begin{theorem}}
\newcommand{\et}{\end{theorem}}
\newcommand{\bc}{\begin{cor}}
\newcommand{\ec}{\end{cor}}
\newcommand{\bl}{\begin{lemma}}
\newcommand{\el}{\end{lemma}}
\newcommand{\bx}{\begin{example}}
\newcommand{\ex}{\end{example}}
\newcommand{\bea}{\begin{eqnarray}}
\newcommand{\eea}{\end{eqnarray}}
\newcommand{\ben}{\begin{eqnarray*}}
\newcommand{\een}{\end{eqnarray*}}

\newcommand{\be}{\begin{equation}}
\newcommand{\ee}{\end{equation}}

\usepackage{color}
\definecolor{cit}{rgb}{1,0,0.2}
\definecolor{link}{rgb}{0.8,0,0}
\usepackage[linktocpage=true,colorlinks,linkcolor=link,citecolor=cit,pagebackref=true]{hyperref}
\begin{document}
\begin{center}
{\bf\Large The Optimal Dynamic Reinsurance Strategies in Multidimensional Portfolio}\\
{Khaled Masoumifard}, {Mohammad Zokaei}\\
\today
\end{center}
\begin{abstract}
The present paper addresses the issue of choosing an optimal dynamic reinsurance policy, which is state-dependent, for an insurance company that operates under multiple insurance business lines. For each line,  the Cramer-Landberg model is adopted for the risk process and one of the contracts such as
 Proportional reinsurance, excess-of-loss reinsurance (XL) and limited XL reinsurance (LXL)  is intended for transferring a portion of the risk to reinsurance.
 In the optimization method used in this paper, the survival function is maximized relative to the dynamic reinsurance strategies. 
 The optimal survival function is characterized as the unique nondecreasing viscosity solution of the associated Hamilton-Jacobi-‌‌‌‌‌‌Bellman equation (HJB) equation with limit one at infinity.  
 The finite difference method (FDM) has been utilized for the numerical solution of the optimal survival function and optimal dynamic reinsurance strategies and the proof for the convergence of the numerical solution to the survival probability function is provided.
  The findings of this article provide insights for the insurance companies as such that based upon the lines in which they are operating, they can choose a vector of the optimal dynamic reinsurance strategies and consequently transfer  some part of their risks to several reinsurers. Using numerical examples, the significance of the elicited results in reducing the probability of ruin is demonstrated in comparison with the previous findings.
\end{abstract}

\NI {\bf Keywords:} Cramer-lundberg process; Optimal reinsurance; Hamilton-Jacobi-Bellman equation; Viscosity solution; Dynamic programming principle.
\section{Introduction}
An effective way for an insurance firm to manage its risk is to buy reinsurance. According to the reinsurance contract, some parts of the claim are shared with the reinsurer, and against, the insurance firm pays a part of its income premium to the reinsurance.
Determination of the optimal reinsurance contracts has been discussed extensively in the literature.  Dynamic proportional reinsurance in the classical risk model for the minimization of the ruin probability was first studied by  \cite{schmidli2001optimal}. \cite{hipp2003optimal} utilized the concept of dynamic excess of loss  reinsurance to extend the Schimidli approach.  \cite{schmidli2002minimizing} have suggested the optimal investment and reinsurance strategies to minimise the ruin probability and have concluded that the investment and reinsurance decrease the ruin
probability for larger initial surplus under the Pareto claim sizes. In this direction, \cite{taksar2003optimal} and \cite{schmidli2004cramer} developed the above
approach  to the diffusion model. Subsequently,
\cite{irgens2004optimal} discussed the maximization of the expected utility of the asset of an insurance
company under the reinsurance and investment constraints in a diffusion classical risk model. A general presentation on ruin probability minimization by means of reinsurance in a classical and diffusion risk models can be found in \cite{schmidli2007stochastic}. Some  additional results with a focus on non-proportional reinsurance contracts are given in \cite{hipp2010optimal}. Recently,  \cite{cani2017optimal} studied a dynamic reinsurance problem obtained from an economical valuation criterion in risk theory  introduced by \cite{hojgaard1998optimal,hojgaard1998optimal1}.

In this paper, we assume that the insurance company produces multiple types of coverage, where customers may purchase different types of insurance policies (such as fire, health, vehicle, etc.). Due to the different risk processes in different lines, it is reasonable that the insurance companies use several reinsurances to share their risk.  For instance, it is possible for an insurance company to  purchase a proportional reinsurance in one line and an excess-of-loss reinsurance in another line or buy one type of excess-of-loss in one line and a different type of excess-of-loss insurance in another line. The survival function, in this paper, is considered as the objective function, and the vector of reinsurance strategies is obtained in such a way that the objective function is maximized; therefore, the results presented in \cite{azcue2014stochastic}, which use a dynamic reinsurance strategy for transferring risk to reinsurers, have been generalized in such a way that the insurer uses the vector of dynamic reinsurance strategies to transfer risk to several reinsurers.

In the second section of the paper, a brief introduction of the model with the presence of the reinsurer and the statement of the problem are provided. In the third section, the main results and in the fourth section numerical examples are presented.
\section{Model formulation}
In the classical Cramer-Lundberg process, the reserve $X_t$ of an insurance company can be described by
\begin{eqnarray}
X_t=x+pt-\sum_{i=1}^{N_t}U_i
\end{eqnarray}
where $N_t$ is a Poisson process with claim arrival intensity $\beta >0$ and the claims size $U_i$ are i.i.d random variables with distribution $F$. The premium rate $p$ is calculated using the expected value principle with relative safety loading $\eta>0 $; that is, $p=(1+\eta)\beta\mu$.the limitation of this model is the assumption that insurers produce only one type of insurance, but in practice, most insurers produce different types of coverage. (e.g. automobile insurance,  fire insurance, workers' compensation insurance, etc.). 
The idea for modeling the surplus process for a company that produces multiple types of coverage  is as follows:
consider the process  $X_t=(X_t^1,\cdots,X_t^n)$  defined as;
\begin{eqnarray}
X_t^{(k)}=p_kt-\sum_{i=1}^{N_t^{(k)}}U_i^{(k)},\qquad \qquad \qquad k=1,\cdots,n
\end{eqnarray}
where $N_t^{(k)}$ is a Poisson process with claim arrival intensity $\beta_k >0$ and the claims size $U_i^{(k)}$'s are i.i.d random variables with distribution $F_k$. Let the risk process of the $k$th line of insurance company be modelled by $X_t^{(k)}$. The premium rate $p_k$ is calculated using the expected value principle with relative safety loading $\eta>0 $; that is, $p_k=(1+\eta)\beta_k\mu_k$. Given an initial surplus $x$, the surplus $Y_t$ of the insurance company at time $t$ can be written as
$
Y_t=x+\sum_{k=1}^{n}X_t^{(k)}
$
and if $X_t^{(1)},\cdots,X_t^{(n)}$ are independent random variables, then $Y_t$ has a compound Poisson distribution, that is,
\begin{eqnarray}
Y_t=x+(\sum_{k=1}^{n}p_k)t-\sum_{i=1}^{N_t}U_i,\qquad \qquad \qquad k=1,\cdots,n
\end{eqnarray}
where $N_t$ is a Poisson process with claim arrival intensity $\beta=\sum_{i=1}^{n}\beta_i$ and the claims size $U_i$'s are i.i.d random variables with distribution $\sum_{i=1}^{n}\frac{\beta_i}{\sum_{i=1}^{n} \beta_i}F_i$.
Let 
$(\Omega_k, \Sigma_k, ({{\cal{F}}_k}_t)_{t \ge 0}, P_k)$ be
the filtered probability space corresponds to line $k$, then, we can describe  filtered probability space model by 
\begin{align}\label{2.3}
(\boldsymbol{\Omega}, \boldsymbol{\Sigma}, (\boldsymbol{{\cal{F}}}_t)_{t \ge 0}, \boldsymbol{P})=(\Omega_1, \Sigma_1, ({{\cal{F}}_1}_t)_{t \ge 0}, P_1)\times \cdots \times(\Omega_n, \Sigma_n, ({{\cal{F}}_n}_t)_{t \ge 0}, P_n).\end{align}

Reinsurance can be an effective way to manage risk by transferring risk from an
insurer to a second insurer (referred to as
the reinsurer).
A reinsurance contract is an agreement between an insurer and a reinsurer under which, claims that arise are shared between the insurer and reinsurer.

Let a Borel measurable function $R:[0,\infty)\longrightarrow[0,\infty)$, called retained loss function, describing the part of the claim that the company pays and satisfies $0\le R(\alpha)\le \alpha$. The reinsurance company covers $\alpha - R(\alpha)$), where the size of the claim is $\alpha$.  Now assume that  in order to reduce the risk exposure of the portfolio, the insurer has the possibility to take reinsurances in a dynamic way for some insurance lines, each of these reinsurances is indexed by $\{1,\cdots,n\}$.  We denote by $\boldsymbol{\mathcal{R}}$ the vector $(\mathcal{R}_1,\cdots,\mathcal{R}_n)$, in which $\mathcal{R}_i$ is the family of retained loss functions associated to the reinsurance policy in $i$'th line, and denote by 
$\Pi_{x}^{\boldsymbol{R}}$ the set of all control strategies with initial surplus $x \ge 0$. 
So, the reinsurances control strategy is a collection 
 $\boldsymbol{R}=(\boldsymbol{R}_t)_{t \ge 0}=({R_1}_t,\cdots, {R_n}_t)_{t \ge 0}$
of the {\Large }vector functions
 $\boldsymbol{R}_t:\boldsymbol{\Omega}\rightarrow \boldsymbol{\mathcal{R}}$  for any $t\ge 0$. 
 
 Well-known reinsurance types are:
 \begin{itemize}
\item[(1)] Proportional reinsurance with $R_P(\alpha)=b\alpha$,
 	\item[(2)] Excess-of- loss reinsurance (XL) with $R_{XL}(\alpha)=\min(\alpha,M)$, $0\le M \le \infty$, 
 	\item[(3)]Limited XL reinsurance (LXL) with $R_{LXL}(\alpha)=\min\{(\alpha,M)\}+(\alpha-M-L)^+$, $0\le M,L \le \infty$.
 \end{itemize}
 In this paper, we assume that the reinsurance calculates its premium using the expected value principle with reinsurance safety loading factor $\eta_1 \ge \eta >0$,
 $$q_R=(1+\eta_1)\beta E(U_i-R(U_i))$$
 and so $p_R=\sum{k} p_k-q_R$. 
  Now, the surplus process in the presence of  reinsurance strategies can be written as
 \begin{eqnarray}\label{2.3}
 X_t^{\boldsymbol{R}}=&x+\sum_{k=1}^{n}\left(\int_{0}^{t}p_{{R_k}_s}ds-\sum_{i=1}^{N_t^{(k)}}{R_k}_{\tau_i}(U_i^{(k)})\right)\nonumber\\
 =& x+\int_{0}^{t}\sum_{k=1}^{n} p_{{R_k}_s}ds-\sum_{k=1}^{n}\sum_{i=1}^{N_t^{(k)}}{R_k}_{\tau_i}(U_i^{(k)}).
 \end{eqnarray}
Without losing the generality, we can write
\begin{eqnarray}\label{2.3'}
X_t^{\boldsymbol{R}}=x+\int_{0}^{t}\sum_{k=1}^{n} p_{{R_k}_s}ds-\sum_{i=1}^{N_t}Z_i
\end{eqnarray}
where  $Z_i$'s are i.i.d random variables with distribution $F_{\boldsymbol{R}}=\sum_{i=1}^{n}\frac{\beta_i}{\sum_{i=1}^{n} \beta_i}F_{R_i}$.
The time of ruin for this process is defined by 
\begin{align}
\label{2.4}
\tau^{\boldsymbol{R}}=\inf\left\{t \ge 0: X_t^{\boldsymbol{R}}<0 \right\}.
\end{align}
In this paper , we aim to identify the reinsurance strategies  $\boldsymbol{R}=(\boldsymbol{R}_t)_{t \ge 0}=({R_1}_t,\cdots, {R_n}_t)_{t \ge 0}$ that maximize survival probability
$ \delta^{\boldsymbol{{R}}}(x)=P(\tau^{\boldsymbol{{R}}}=\infty|X_0^{\boldsymbol{{R}}}=x)$
, in other words, we are looking for the optimal survival function
\begin{eqnarray}\label{delta}
\delta(x)=\underset{\boldsymbol{R}\in \Pi_x^{\boldsymbol{R}}}{\sup}\,\,\delta^{\boldsymbol{R}}(x).
\end{eqnarray}
From now forward, where we use $\delta(.)$, we take it to mean the optimal survival function.
It is easy to show, similar to section 2.1.1 of \cite{azcue2014stochastic}, that the HJB equation of this problem can be written as
\begin{eqnarray}\label{hjbeq}
\underset{\boldsymbol{\mathcal{R}}}{sup}\,\,{\cal{L}}_{\boldsymbol{R}}(\delta)(x)=0
\end{eqnarray}
where 
\begin{align}\label{Lfunction}
{\cal{L}}_{\boldsymbol{R}}(\delta)(x)=p_{\boldsymbol{R}}\delta'(x)-(\sum_{i=1}^{n}\beta_i)\delta(x)+(\sum_{i=1}^{n}\beta_i)\int_{0}^{x}\delta(x-\alpha)dF_{\boldsymbol{R}}(\alpha).
\end{align}
\section{Main  results}
The dynamic programming method is a cogent means
to scrutinize the stochastic control problems through  the HJB  equation.
In (\ref{hjbeq}),  we have obtained the associated equation to the
value function (\ref{delta}). Nonetheless, in the classical approach, this method is adopted only when it is assumed \textit{a priori} 
that optimal value functions are smooth enough. In general,  the optimal value function is not expected  to be smooth enough to satisfy these equations in the classical sense. 
 These call for a felt need for a week notation of solution of the HJB equation: the theory of viscosity solutions. Let us define this notion(see \cite{azcue2014stochastic}).
\begin{defn}
	We say that a  function $\underline{u}:[0,\infty)\rightarrow \mathbb{R} $ is a viscosity subsolution of (\ref{hjbeq}) at $x \in (0,\infty)$ if it is locally Lipschitz and any continuously differentiable function $\varphi:(0,\infty)\rightarrow \mathbb{R} $ (called test function) with $\varphi(x)=\underline{u}(x)$ such that $\underline{u}-\varphi$ reaches the maximum at $x$ and satisfies 
	$\underset{\boldsymbol{\mathcal{R}}}{sup}\,\,{\cal{L}}_{\boldsymbol{{R}}}(\varphi)(x)\ge 0$. 	We say that a continuous function $\bar{u}:[0,\infty)\rightarrow \mathbb{R}$ is a viscosity supersolution of (\ref{hjbeq}) at $x \in (0,\infty)$ if it is locally Lipschitz  any continuously differentiable function $\phi:(0,\infty)\rightarrow \mathbb{R}$ (called test function)  with $\phi(x)=\bar{u}(x)$ such that $\bar{u}-\phi$ reaches the maximum at $x$ and satisfies 
	$\underset{\boldsymbol{\mathcal{R}}}{sup}\,\,{\cal{L}}_{{\boldsymbol{R}}}(\phi)(x)\le 0$.
	Finally, we say that a continuous function $u:[0,\infty)\rightarrow \mathbb{R}$ is a viscosity solution of (\ref{hjbeq}) if it is both a viscosity subsolution and a viscosity supersolution at any $x \in (0,\infty)$.
\end{defn}
The methods which will  be used later are gleaned from  \cite{azcue2014stochastic} and \cite{nozadi2014optimal}.
\subsection{ Viscosity Solution}
In chapter 3 of \cite{azcue2014stochastic}, there is an example that the survival probability function in problem with reinsurance can be non-differentiable. Hence, in general, we cannot expect for $\delta$ to have the smoothness properties needed to regard it as a solution (in the classical sense) for the corresponding HJB equation (\ref{hjbeq}). We prove instead that
$\delta$ is a viscosity solution of the corresponding HJB equation. 
Before stating the main results,  the following lemmas are needed.

\begin{lemma}\label{lemma2.1}
Consider an arbitrary admissible strategy $\boldsymbol{R}=(\boldsymbol{R}_t)_{t\ge 0}\in \Pi_{x}^{\boldsymbol{R}}$ and set  $x \ge 0$. Then $\delta^{\boldsymbol{R}}\big(X_{t \wedge \tau^{\boldsymbol{R}}\wedge t}^{\boldsymbol{R}}\big)$ is a martingale.
\end{lemma}
\begin{lemma}\label{S-lipschitz}
	 $0\le \delta(x) \le 1$ for all $x \ge 0$, $\lim_{x\rightarrow \infty} \delta(x)=1$, $\delta$ is increasing, and it is Lipschitz with Lipschitz constant $K=\beta/(\sup_{\boldsymbol{R}\in \boldsymbol{\mathcal{R}}} \sum_{i=1}^{n} p_{R_i})$.
\end{lemma}
\begin{lemma}\label{Sprop2.13}
	Let a nonnegative continuously differentiable function $u:\mathbb{R}_+\rightarrow \mathbb{R}$. Then for any finite stopping time $\tau\le \tau^{\boldsymbol{R}}$, the following equality is true,
	$$u(X^{\boldsymbol{R}}(t))-u(x)=\int_{0}^{\tau}{{\cal{L}}}_{\boldsymbol{\mathcal{R}}}(u)(X_{s^-}^{\boldsymbol{R}})ds+M_{\tau}$$
	where $X_{\tau}^{\boldsymbol{R}}$, $\tau^{\boldsymbol{R}}$ and $ {{\cal{L}}}_{\boldsymbol{\mathcal{R}}}$ are defined in  (\ref{2.3}), (\ref{2.4}) and (\ref{Lfunction}) and $M_{\tau}$ is a martingal with zero expectation.
\end{lemma}
 \begin{lemma}\label{lemma3.2.'}
 	Let the vector $\boldsymbol{\mathcal{R}}=(\mathcal{R}_1,\cdots,\mathcal{R}_n)$, where $\mathcal{R}_i$ is one of the reinsurance families $\mathcal{R}_P$, $\mathcal{R}_{XL}$ and $\mathcal{R}_{LXL}$. 
 	If $u$ is a nonnegative and a twice continuously
 	differentiable function defined in $\boldmath{R}_+$ (extended as $u(x)=0$ for $x < 0$), then $H=\underset{\boldsymbol{R} \in \boldsymbol{\mathcal{R}}}{sup} {\cal{L}}_{\boldsymbol{R}}(u)(.)$ is upper semicontinuous. Moreover,  for any $A > 0$ and $h \in (0,1)$, there exists $K_A > 0$ such that
 	$$|H(y)-H(x_0)|\le K_{A}(y-x_0+F(y)-F(x_0)), \qquad \qquad 0\le x_0,\,\, y\le A.$$
 	 \end{lemma}
Now the following theorem is obtained.
\begin{theorem}\label{prop1}
	Let the vector $\boldsymbol{\mathcal{R}}=(\mathcal{R}_1,\cdots,\mathcal{R}_n)$, where $\mathcal{R}_i$ is one of the reinsurance families $\mathcal{R}_P$, $\mathcal{R}_{XL}$ and $\mathcal{R}_{LXL}$. Then, the function $\delta$  is a viscosity solution of (\ref{hjbeq}).
\end{theorem}
{\bf  Proof }\,\,
The proof of this theorem is quite similar to  the Proposition 3.2 of \cite{azcue2014stochastic}. Indeed, applying the Lemmas \ref{lemma2.1}, \ref{S-lipschitz}, \ref{Sprop2.13} and \ref{lemma3.2.'} and repeating the proof of the Proposition 3.2 of \cite{azcue2014stochastic} provide the proof. 
$\square$
\subsection{Characterization}
In this section, it will be proved that the probability function defined in (\ref {delta}) is the unique viscosity solution of the HJB   (\ref {hjbeq}) with limit one at infinity. In order to prove the uniqueness result, we use the following three lemmas.
 \begin{lemma}\label{lemm4.2}
Suppose that 
	 $\bar{u}$
is a non-decreasing Lipschitz viscosity supersolution of
	(\ref{hjbeq}) 
	such that
	$$\lim_{x\rightarrow \infty}\bar{u}(x)=1.$$ 
	A sequence of positive functions $\bar{u}_n:\mathbb{R}_+\rightarrow \mathbb{R}$ are found such that
	\begin{itemize}
		\item[(1)] $\bar{u}_n$
is infinitely continuously differentiable and 
		$\bar{u}_n'\le K$, where $K$ is the Lipschitz constant of
		$\bar{u}$.
		\item[(2)] $\lim_{x\rightarrow \infty} \bar{u}_n(x)=1$
		\item[(3)] $\bar{u}_n \searrow \bar{u}$ uniformly in  
		$\mathbb{R}_+$ 
and
		$\bar{u}_n'(x)$
converges to 
		$\bar{u}'(x)$ a.e.
		\item[(4)] 
There is a sequence
		$c_n >0$
with
		$\lim_{n\rightarrow \infty} c_n = 0$
such that
		$$\sup_{\boldsymbol{R}\in \boldsymbol{\mathcal{R}}, x \ge 0} 	{\cal{L}}_{\boldsymbol{R}}(\bar{u}_n)(x)\le c_n.$$
		
	\end{itemize}
\end{lemma}
\begin{lemma}\label{lemm4.3}
Suppose that 
	 $\underline{u}$
is a non-decreasing Lipschitz viscosity supsolution for
	(\ref{hjbeq}), 
	such that
	$$\lim_{x\rightarrow \infty}\underline{u}(x)=1.$$
	A sequence of positive functions $\underline{u}_n:\mathbb{R}_+\rightarrow \mathbb{R}$ are found such that
	\begin{itemize}
		\item[(1)] $\underline{u}_n$
is infinitely continuously differentiable and 
		$\underline{u}_n'\le K$, where $K$ is the Lipschitz constant of
		$\underline{u}$.
		\item[(2)] $\lim_{x\rightarrow \infty} \underline{u}_n(x)=1$
		\item[(3)] $\underline{u}_n \nearrow \underline{u}$ uniformly in  
		$\mathbb{R}_+$ 
and
		$\underline{u}_n'(x)$
converges to 
		$\underline{u}'(x)$ a.e.
		\item[(4)] 
There is a sequence
		$c_n >0$
with
		$\lim_{n\rightarrow \infty} c_n = 0$
such that
		$$\sup_{\boldsymbol{R}\in \boldsymbol{\mathcal{R}}, x \ge 0} 	{\cal{L}}_{\boldsymbol{R}}(\underline{u}_n)(x)\ge -c_n.$$
		
	\end{itemize}
\end{lemma}
\begin{lemma}\label{Lemma4.4}
	Consider twice continuously differentiable function $\bar{u}$. Then, for $\epsilon >0,$ there exists a stationary reinsurance strategy $\boldsymbol{\rho}=(\rho_1,\cdots,\rho_n)$  satisfying the following,
	$$\sup_{\boldsymbol{R} \in \boldsymbol{\mathcal{R}}} {\cal{L}}_{\boldsymbol{R}}(x)- {\cal{L}}_{\boldsymbol{\rho}^x}(x)<\epsilon,\qquad x \ge 0.$$
\end{lemma}
\begin{prop}\label{Smal_vis_Large_visc}
$\delta(.)$  is both the smallest viscosity supersolution and the largest viscosity sub-solution of HJB (\ref{hjbeq}), with limit one at infinity.
\end{prop}
{\bf  Proof }\,\,
To prove this theorem, it suffices to repeat the proof of Theorem 4.3 of \cite{azcue2014stochastic} by replacing   $\boldsymbol{R}=(\boldsymbol{R}_t)_{t \ge 0}=({R_1}_t,\cdots, {R_n}_t)_{t \ge 0}$ with  $\bar{R}$.
At first, it is shown that $\delta(.)$  is smaller or equal to that any supersolution. Assume  that 
$\bar{u} $ is a non-decreasing viscosity supersolution satisfying in 
(\ref{hjbeq}), with 
$\lim_{x\rightarrow \infty} \bar{u}(x)=1$.
Consider an arbitrary admissible strategy $\boldsymbol{R}=(\boldsymbol{R}_t)_{t\ge 0}\in \Pi_{x}^{\boldsymbol{R}}$ and   
$ x \ge 0$. Denote by $X_t$, the controlled risk process with initial surplus  
$x$ corresponding to $\boldsymbol{R}$ and let $\tau$  be its ruin time. For any  $M>x$, define the following stopping time
$$\tau_M=\inf\{t\ge 0:X_t^{\boldsymbol{R}} \ge M\}.$$
Consider the function 
$\bar{u}_n$ defined in Lemma  
\ref{lemm4.2},  and set $\bar{u}_n=0$ in 
$(-\infty, 0)$. From Lemma
\ref{Sprop2.13} and part (4) of Lemma  
\ref{lemm4.2}, it follow that:
\begin{align*}
\bar{u}_n(X_{\tau_M\wedge \tau \wedge t})-\bar{u}_n(x)&=\int_{0}^{\tau_M\wedge \tau \wedge t}{\cal{L}}_{\boldsymbol{R}_s}( \bar{u}_n)(X_{s^-})ds+M_{\tau_M\wedge \tau \wedge t}\\
&\le c_n({\tau_M\wedge \tau \wedge t})+M_{\tau_M\wedge \tau \wedge t}
\end{align*} 
wherein 
$M_t$ 
is martingale with zero mean. Thus, the following inequality is obtained
\begin{align*}
E_x( \bar{u}_n(X_{\tau_M\wedge \tau \wedge t}))-\bar{u}_n(x)\le c_n t. 
\end{align*}
Now, taking limit from both sides of the above relation when
$n\rightarrow \infty$, for fix $t$, the next inequality is established:
$$\limsup_{n\rightarrow \infty} E_x( \bar{u}_n(X_{\tau_M\wedge \tau \wedge t}))-\bar{u}_n(x)\le 0$$
and so, since  
$\lim_{n\rightarrow \infty}\bar{u}_n(x)=\bar{u}(x)$ 
and
$\bar{u}(x)\le \bar{u}_n(x)$, 
$ E_x( \bar{u}(X_{\tau_M\wedge \tau \wedge t}))-\bar{u}(x)\le 0.$ 
Therefore, when $t\rightarrow \infty$,  the following inequality is obtained
$$ E_x( \bar{u}(X_{\tau_M\wedge \tau \wedge t}))-\bar{u}(x)=\bar{u}(M)P(\tau_M<\tau)-\bar{u}(x)\le 0.$$ 
Now, taking 
$M\rightarrow \infty$, and using the fact that, 
$\lim_{M\rightarrow \infty}\bar{u}(M)=1$ 
and 
$\lim_{M\rightarrow \infty}P(\tau_M < \tau)=\delta^{\boldsymbol{R}}(x)$
it follows that 
$$\delta^{\boldsymbol{R}}(x) \le \bar{u}(x)\qquad for \,\, all\,\,\boldsymbol{R}\in \Pi_{x}^{\boldsymbol{R}},$$
which in turn results in $\delta(x) \le \bar{u}(x).$\\
We shall now prove that $\delta(.)$  is greater or equal than any subsolution. Assume that 
$\underline{u}$
is a non-decreasing subsolution in 
(\ref{hjbeq}) with
$\lim_{x\rightarrow \infty} \underline{u}(x)=1$. Consider the functions  
$\underline{u}_n$ defined in Lemma 
\ref{lemm4.3}, and set $\underline{u}_n=0$  in 
$(-\infty,0)$. Based on Lemma  
\ref{Lemma4.4}, for any 
$y\ge 0$ and 
$n\ge 1$, there exists a stationary reinsurance control  
$\boldsymbol{\rho}_n$ such that, 
$$ \sup_{\boldsymbol{R}\in \boldsymbol{\mathbb{R}}}{\cal{L}}_{\boldsymbol{R}}(\underline{u})(y)-{\cal{L}}_{\boldsymbol{\rho}_n^y}(\underline{u}_n)(y)\le \frac{1}{n}.$$
Consider the controlled process  
$(X_t^n)_{t\ge 0}$   with initial surplus  $x$ and admissible reinsurance strategy  
$\boldsymbol{R}^n=(\boldsymbol{R}_t^n)_{t\ge 0}=(\boldsymbol{\rho}_n^{X_{t^-}})_{t\ge 0}$ and denote by $\tau^n$, the corresponding ruin time. For each  
$M>0$, define the following stopping time:
$$\tau_M^n=\inf\{t\ge 0:X_t^n\ge M\}.$$ 
From Lemma 
\ref{Sprop2.13} 
and Lemma 
\ref{lemm4.3}, it follows that for each $n$, 
\begin{align*}
\underline{u}_n(X_{\tau_M^n\wedge \tau^n \wedge t}^n)-\underline{u}_n(x)&=\int_{0}^{\tau_M^n \wedge \tau^n \wedge t}{\cal{L}}_{\boldsymbol{R}_s^n}( \underline{u}_n)(X_{s^-}^n)ds+M_{\tau_M^n\wedge \tau^n \wedge t}\\
&\ge \int_{0}^{\tau_M^n \wedge \tau^n \wedge t}\sup_{\boldsymbol{R}\in \boldsymbol{\mathbb{R}}}\left({\cal{L}}_{\boldsymbol{R}}( \underline{u}_n)(X_{s^-}^n)-\frac{1}{n}\right) ds+M_{\tau_M^n\wedge \tau^n \wedge t}\\
&\ge \left(-c_n-\frac{1}{n}\right)({\tau_M^n\wedge \tau^n \wedge t})+M_{\tau_M^n\wedge \tau^n \wedge t}, 
\end{align*} 
wherein 
$M_t$  is a martingale with mean zero. So, the following inequality is obtain by
\begin{align*}
E_x\left(\underline{u}_n(X_{\tau_M^n\wedge \tau^n \wedge t}^n)-\underline{u}_n(x)\right)\ge \left(-c_n t-\frac{1}{n}\right)E_x\left(\tau_M^n\wedge \tau^n \wedge t\right). 
\end{align*}
Now, taking limit from both sides of the above relation when
$n\rightarrow \infty$, for a fix $t$, the next inequality is established as,
$$\limsup_{n\rightarrow \infty} E_x\left(\underline{u}_n(X_{\tau_M^n\wedge \tau^n \wedge t}^n)\right)\ge \underline{u}(x)$$
and so, since  
$\underline{u}(x)\ge \underline{u}_n(x)$, 
$$\limsup_{n\rightarrow \infty} E_x\left(\underline{u}(X_{\tau_M^n\wedge \tau^n \wedge t}^n)\right)\ge \underline{u}(x).$$ 
For  
$\epsilon >0$, take 
$n_0$ 
large enough so that
$$E_x\left(\underline{u}(X_{\tau_M^{n_0}\wedge \tau^{n_0} \wedge t}^{n_0})\right)\ge \underline{u}(x)-\epsilon.$$
Thus, 
\begin{align*}
\underline{u}(x)-\epsilon &\le E_x\left(\underline{u}(X_{\tau_M^{n_0}\wedge \tau^{n_0} \wedge t}^{n_0})\right)\\
&=\underline{u}(M)P(\tau_M^{n_0}\wedge \tau^{n_0} \wedge t)+E_{x}\left(\underline{u}(X_{\tau_M^{n_0}\wedge \tau^{n_0} \wedge t}^{n_0})I_{t\le \tau_M^{n_0}\wedge \tau^{n_0} }\right). \end{align*}
Note that 
$P(\tau_M^{n_0}\wedge \tau^{n_0} \wedge t)$ is non-decreasing in $t$,   
$0\le \underline{u}\le 1$,
and
$\lim_{t\rightarrow \infty}P(\tau_M^{n_0}\wedge \tau^{n_0} > t)=0,$
 which result in 
$$\underline{u}(x)-\epsilon\le \liminf_{t\rightarrow \infty} E_x\left(\underline{u}(X_{\tau_M^{n_0}\wedge \tau^{n_0} \wedge t}^{n_0})\right)
=\underline{u}(M)P(\tau_M^{n_0}< \tau^{n_0} ).$$
Now, taking  
$M\rightarrow \infty$,  
the following result is established:
$$\underline{u}(x)-\epsilon\le \delta^{\boldsymbol{R}^{n_0}}\le \delta(x)$$ and hence 
$\underline{u}(x)\le \delta(x)$.$\square$

In the next theorem, the optimal survival function is characterized.
the theorem  provided below is the cospicuously obvious result of Theorem \ref*{prop1} and Proposition \ref{Smal_vis_Large_visc} .
\begin{theorem}\label{uniquesol}
	$\delta(.)$ is the unique nondecreasing viscosity solution of (\ref{hjbeq}) with limit one at infinity. 
\end{theorem}
We summarize these results in the following corollary.
\begin{cor}
(Verification) If the survival probability function of a vector of reinsurance admissible strategies is a viscosity solution of the HJB equation (\ref {hjbeq}) with limit one at infinity, then the vector of reinsurance admissible strategies and its survival probability function are optimal.
\end{cor}

\subsection{Numerical solution }
Let us define the function
 $$\hat{\delta}(x)=\inf_{\boldsymbol{R}\in \boldsymbol{\mathcal{R}}}\frac{(\sum_{i=1}^{n}\beta_i)\delta(x)-(\sum_{i=1}^{n}\beta_i)\int_{0}^{x}\delta(x-\alpha)dF_{\boldsymbol{R}}(\alpha)}{p_{\boldsymbol{R}}}.$$
 Similarly to the Lemma 5.6 of \cite{azcue2014stochastic}, we can show that $\hat{\delta}(x)$  is well defined, nonnegative, Borel measurable, $\hat{\delta}=\delta'$ a.e., and 
 $$\sup_{\boldsymbol{R}\in \boldsymbol{\mathcal{R}}} \big\{p_{\boldsymbol{{R}}}\hat{\delta}(x)-(\sum_{i=1}^{n}\beta_i)\delta(x)+(\sum_{i=1}^{n}\beta_i)\int_{0}^{x}\delta(x-\alpha)dF_{\boldsymbol{{R}}}(\alpha) \big\}=0.$$
In this section, using FDM, we try to solve numerically the problem of optimal reinsurance and optimal survival function. The utilization of the FDM is prevalent for solving the HJB equations in stochastic control problem and these numerical solutions are usually converges to the viscosity solution (see \cite{crandall1992user}, \cite{fleming2006controlled}, \cite{pham2009continuous} and \cite{nozadi2014optimal}).
 A numerical solution for
$
sup_{\boldsymbol{{R}} \in \boldsymbol{\mathcal{R}}}\,\,{\cal{L}}_{\boldsymbol{\mathcal{R}}}(f)(x)=0,
$ can be obtained by the use of FDM and adaptation of the boundary condition $\lim_{x\rightarrow \infty} f(x)=1$. Similar to the procedure described in \cite{fleming2006controlled}, section IX.3 or \cite{nozadi2014optimal}, section 3.3, we discretize the state space  with sufficiently small step size $h$  and  define a family of functions $f_h$ in the following procedure: starting with  $$f_h(0)= 1\, \text{and}\, f'_h(0)=\inf_{\boldsymbol{{R}} \in \boldsymbol{\mathcal{R}}}\beta\frac{\
1-p(\boldsymbol{{R}}(Y)=0)}{p_{\boldsymbol{{R}}}}$$  
and for $s=ih,\, i=1,2,\cdots$, we approximate 
  $\int_{0}^{x}f(x-\alpha) dF_R(\alpha)$
by  $$G_{\boldsymbol{R}}(s)=\sum_{\{j\le i\}}f_{h}((i-j)h)P\{(j-1)h<\boldsymbol{R}(Y)\le jh\}.$$
It is easy to show that $G_{\boldsymbol{R}}(s)$ converges to  $\int_{0}^{s}f(x-\alpha) dG_R(\alpha)$ as $h$ tends zero.
Then we define 
  $f'_h(s)$ by
  \begin{align}\label{3.24}
  f'_{h}(s)&=\inf_{{\boldsymbol{R}}}\frac{\beta h(f_{h}(s-h)-G_{{\boldsymbol{R}}}(s))}{{p_{{\boldsymbol{R}}}}}
  \end{align}
and  set $f_h(s)=f_h(s-h)+hf'_h(s)$. 
\begin{lemma}\label{numeric-increasing}
 Let some small step size $ h $ such that $p_{{\boldsymbol{R}}}\ge 0, \,\, i=1,2, \cdots,$  and let $D_h= \{ih, i=0,\,1,\cdots\}$. Then 
 \begin{itemize}
 	\item[(i)]  $f'_h(s) \ge 0$ for all $s\in D_h$,
 	\item[(ii)] for all $k\ge 0$, the following inequalities are true,
 	$$f_h(kh) \le (1-\frac{\beta}{p_{{\boldsymbol{R}}}}h)^{-k}\le e^{\frac{\beta}{p_{{\boldsymbol{R}}}}k h}$$ and $$f'_h(k h)\le\frac{\beta}{p_{{\boldsymbol{R}}}}f_h(kh) .$$
 \end{itemize}
\end{lemma}
In the next Proposition, we use the same argument as in \cite{nozadi2014optimal}, section 3.3, to demonstrate that the function $f_{h}/{f_{h}(\infty)}$ converges to the unique viscosity $\delta$.
\begin{prop}\label{numeric-visco}
In the setting of the above lemma, and define
\begin{align}
f^{*}(s)=\lim_{h\rightarrow 0}\sup_{ih \rightarrow s} f_h(ih),
\end{align}
and 
\begin{align}
f_{*}(s)=\lim_{h\rightarrow 0}\inf_{ih \rightarrow s} f_h(ih).
\end{align}
Then the functions $f^{*}(s)$ and $f_{*}(s)$ are respectively, sub- and super viscosity solution of (\ref{hjbeq}).
\end{prop}
{\bf proof}  Firstly, we show that the  function $f^*$ is locally Lipschitz and a viscosity subsolution of \ref{hjbeq}. Fix $M>0$ and let $0\le x < y\le M$ be arbitrary and take a sequence $h_{n} \rightarrow 0$  as $n \rightarrow \infty$ such that:
 for any number $\epsilon >0$  there exists some number $n_0$ and two sequences $k_{n}^{(1)}\in \mathbb{N}$ and $k_{n}^{(1)}\in \mathbb{N}$ such that  for all $n \ge n_0$ we have   $|f_{h_{n}}(k_{n}^{(2)}h_{n})- f^*(y)|<\epsilon$, $|k_{n}^{(1)}h_{n}-x|< \epsilon$  and $|k_{n}^{(2)}h_{n}-y|< \epsilon$.
 It is easy to see that
 \begin{align*}
 f^*(y)-f^*(x)=& \lim_{n\longrightarrow \infty} f_{h_n}(k_n^{(2)}h_n)-f^*(x)\\
 & \le \lim_{n\longrightarrow \infty} f_{h_n}(k_n^{(2)}h_n)- \lim_{n\longrightarrow \infty}  f_{h_n}(k_n^{(1)}h_n)\\
 & \le \lim_{n\longrightarrow \infty} \big(f_{h_n}(k_n^{(2)}h_n)-f_{h_n}(k_n^{(1)}h_n)\big).
 \end{align*}
 Now, according to  the above inequality and Lemma \ref{numeric-increasing}, the following relation is obtained:
 $$ f^*(y)-f^*(x) \le K(y-x)$$
  where $K$ is a common Lipschitz constant of $f_{h_n}(t)$, $0\le t \le M$.
  
  To show that $f^{*}$ is a viscosity subsolution,  suppose that $w$ is a test function such that $ f^*(x)-w(x)$ has a maximum at $s > 0$.  Therefore
$$ f_{h_{n_0}}(k_{n_0}h_{n_0}+h_{n_0})- w(k_{n_0}h_{n_0}+h_{n_0})\le  f_{h_{n_0}}(k_{n_0}h_{n_0})- w(k_{n_0}h_{n_0}),$$
$$ f_{h_{n_0}}(k_{n_0}h_{n_0}+h_{n_0})- f_{h_{n_0}}(k_{n_0}h_{n_0}) \le w(k_{n_0}h_{n_0}+h_{n_0}) - w(k_{n_0}h_{n_0}).$$
  By definition,
\begin{align*}
0=&\underset{\boldsymbol{\mathcal{R}}}{sup}\,\,\{p_{{\boldsymbol{R}}}f_{h_{n_0}}'(k_{n_0}h_{n_0}+h_{n_0})-\beta f_{h_{n_0}}(k_{n_0}h_{n_0}+h_{n_0})+\beta G_{\boldsymbol{R}}(k_{n_0}h_{n_0}+h_{n_0})\}\\
=&\underset{\boldsymbol{\mathcal{R}}}{sup}\,\,\{p_{{\boldsymbol{R}}}\frac{f_{h_{n_0}}(k_{n_0}h_{n_0}+h_{n_0})- f_{h_{n_0}}(k_{n_0}h_{n_0})}{h_{n_0}}-\beta f_{{h_{n_0}}}(k_{n_0}h_{n_0}+h_{n_0})+\beta G_{\boldsymbol{R}}(k_{n_0}h_{n_0}+h_{n_0})\}\\
\le& \underset{\boldsymbol{\mathcal{R}}}{sup}\,\,\{p_{{\boldsymbol{R}}}\frac{w(k_{n_0}h_{n_0}+h_{n_0}) - w(k_{n_0}h_{n_0})}{h_{n_0}}-\beta f_{h_{n_0}}(k_{n_0}h_{n_0}+h_{n_0})+\beta G_{\boldsymbol{R}}(k_{n_0}h_{n_0}+h_{n_0})\}\\
=& \underset{\boldsymbol{\mathcal{R}}}{sup}\,\,\{p_{{\boldsymbol{R}}}w'(k_{n_0}h_{n_0})-\beta f_{h}(k_{n_0}h_{n_0}+h_{n_0})+\beta G_{\boldsymbol{R}}(k_{n_0}h_{n_0}+h_{n_0})\}.
\end{align*}
Then by Fatou’s lemma,
$$\limsup_{n\Rightarrow \infty} G_{\boldsymbol{R}}(k_{n_0}h_{n_0}+h_{n_0}) \le E_{\boldsymbol{R}}(f^*(s-Y)), \qquad f_{{h_{n_0}}}(k_{n_0}h_{n_0}+h_{n_0})\rightarrow f^*(k_{n_0}h_{n_0}).$$
So,
\begin{align*}
0\le&  \underset{\boldsymbol{\mathcal{R}}}{sup}\,\,\{p_{{\boldsymbol{R}}}w'(k_{n_0}h_{n_0})-\beta f_{{h_{n_0}}}(k_{n_0}h_{n_0}+h_{n_0})(k_{n_0}h_{n_0})+\beta G_{\boldsymbol{R}}(k_{n_0}h_{n_0}+h_{n_0})\}\\
\le & \underset{\boldsymbol{\mathcal{R}}}{sup}\,\,\{p_{{\boldsymbol{R}}}w'(k_{n_0}h_{n_0})-\beta f^*(k_{n_0}h_{n_0})+\beta E_{\boldsymbol{R}}(f^*(x-Y))\}\\
\le & \underset{\boldsymbol{\mathcal{R}}}{sup}\,\,\{p_{{\boldsymbol{R}}}w'(s)-\beta w(s)+\beta E_{\boldsymbol{R}}(w(s-Y))\}.
\end{align*}
Thus, $f ^*$ is a viscosity subsolution. Similarly, $f _*$  is locally Lipschitz and a viscosity supersolution of \ref{hjbeq}.
$\square$
\begin{theorem}
The sequence $f_{h}/{f^*(\infty)}$ converges to the unique viscosity $\delta$.
\end{theorem}
{\bf Proof}
Define $g^*(s)=f^*(s)/{f^*(\infty)}$ and $g^*(s)=f_*(s)/{f_*(\infty)}$.  It is obvious that $\lim_{s \rightarrow \infty}g^*(s)=1$ and by Lemma \ref{numeric-increasing} $g^*(.)$ and $g_*(.)$ are nondecreasing functions. Now by the Proposition \ref{Smal_vis_Large_visc}, $\delta \ge g^*$ and  $\delta \le g_*$ and so  $g^*\le g_*$. Since  $g^*\ge g_*$ by definition, we have convergence. Define $g(s)= g^*(s)$ $(=g_*(s))$. Now, by \ref{numeric-visco}, the function $g(s)$ nondecreasing viscosity solution of (\ref{hjbeq}) with limit one at infinity.  Therefore, by Theorem \ref{uniquesol}, $\delta(s)=g(s)$.  
$\square$
\section{Examples}
\begin{figure}
\centering
\begin{subfigure}{.4\textwidth}
	\caption{}
	\centering
	\includegraphics[width=1\linewidth]{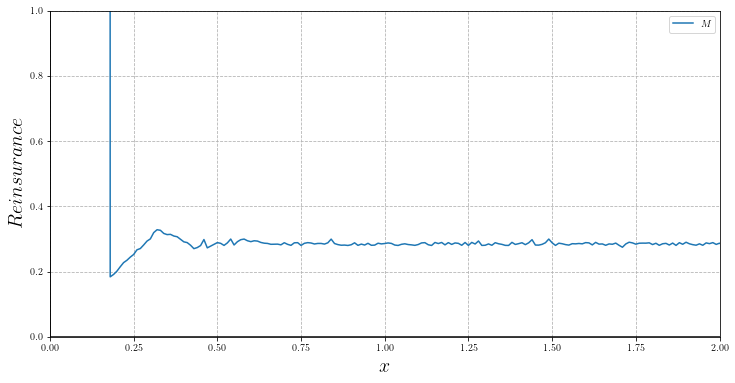}
\end{subfigure}
\begin{subfigure}{.4\textwidth}
	\caption{}
	\centering
	\includegraphics[width=1\linewidth]{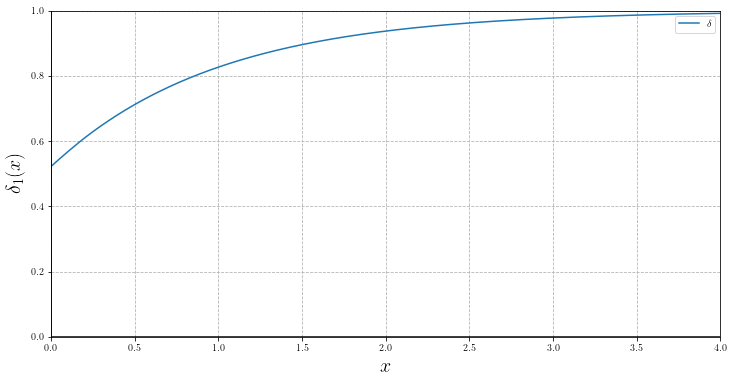}
\end{subfigure}
\caption{Three lines of insurance with one excess of loss reinsurance: (a) Optimal retention levels for excess-of-loss reinsurance (b) Survival probability function with excess-of-loss  reinsurance. }
\label{optimalresultoneXL}
\end{figure}

\begin{figure}
	\centering
		\includegraphics[width=.8\linewidth]{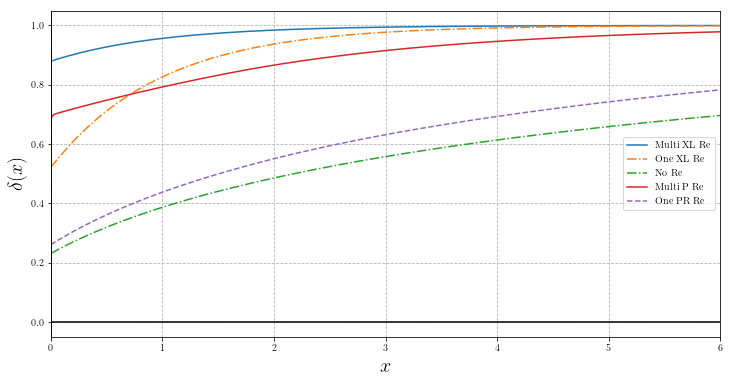}
		\caption{Survival functions}
	\label{optimalresultsurvivalfunctions}
\end{figure}

\begin{figure}
	\centering
	\begin{subfigure}{.4\textwidth}
		\caption{}
		\centering
		\includegraphics[width=0.8\linewidth]{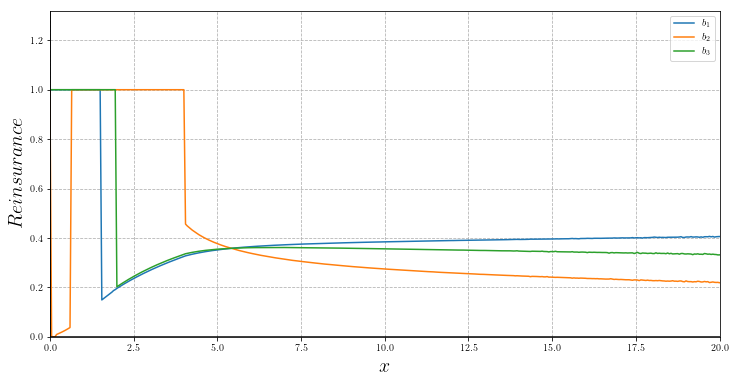}
		\label{ppp}
	\end{subfigure}
	\begin{subfigure}{.4\textwidth}
		\caption{}
		\centering
		\includegraphics[width=0.8\linewidth]{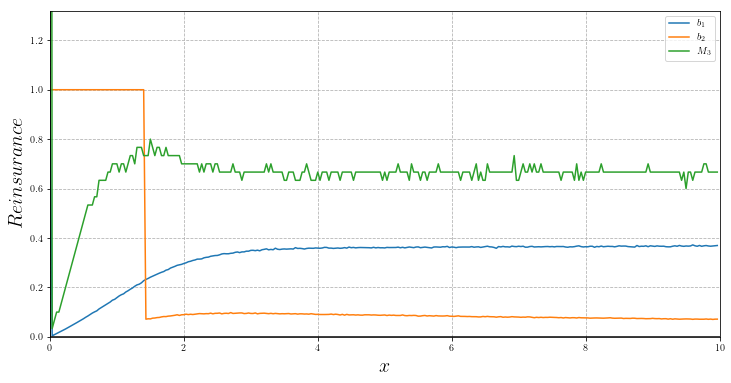}
		\label{ppXL}
	\end{subfigure}
\\

	\begin{subfigure}{.4\textwidth}
		\caption{}
		\centering
		\includegraphics[width=0.8\linewidth]{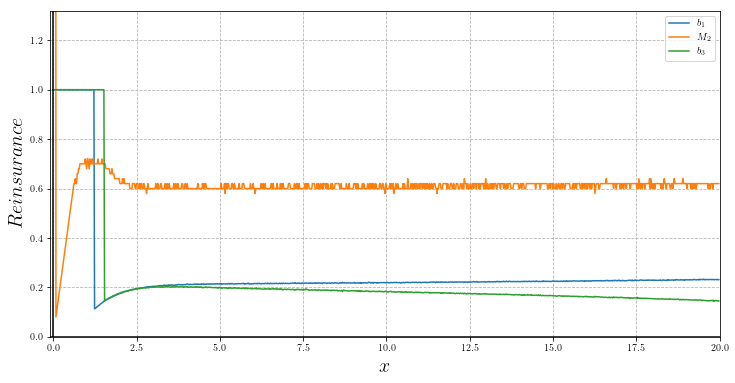}
		\label{pXLp}
	\end{subfigure}
	\begin{subfigure}{.4\textwidth}
		\caption{}
		\centering
		\includegraphics[width=0.8\linewidth]{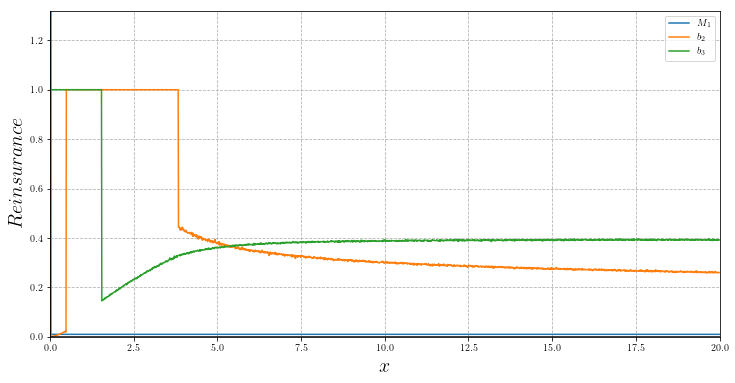}
		\label{XLpp}
	\end{subfigure}
\\
\begin{subfigure}{.4\textwidth}
	\caption{}
	\centering
	\includegraphics[width=0.8\linewidth]{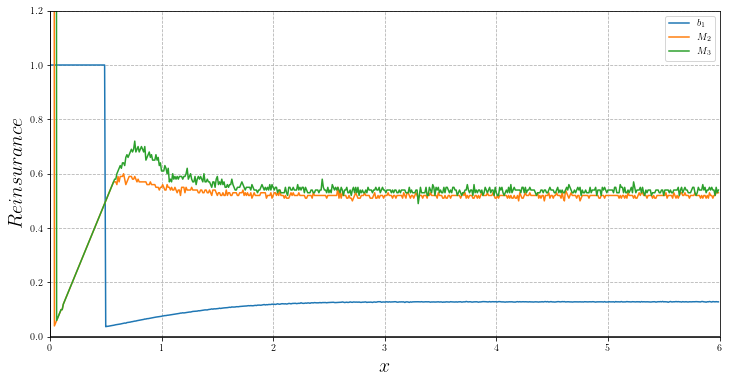}
	\label{pXLXL}
\end{subfigure}
\begin{subfigure}{.4\textwidth}
	\caption{}
	\centering
	\includegraphics[width=0.8\linewidth]{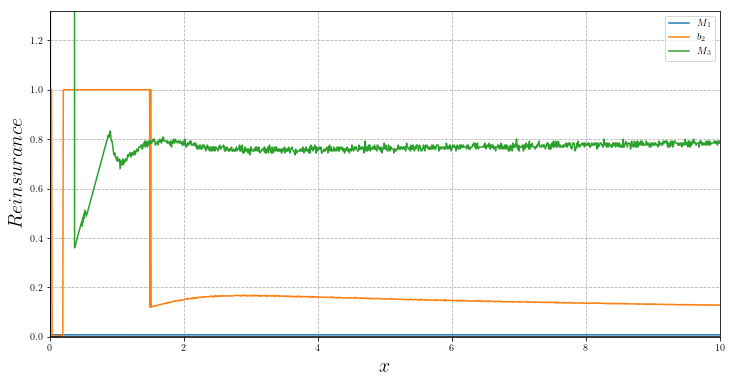}
	\label{XLpXL}
\end{subfigure}
\\
\begin{subfigure}{.4\textwidth}
	\caption{}
	\centering
	\includegraphics[width=0.8\linewidth]{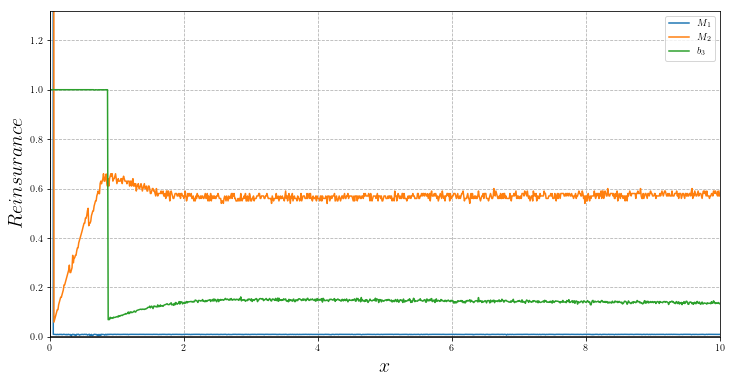}
	\label{XLXLp}
\end{subfigure}
\begin{subfigure}{.4\textwidth}
	\caption{}
	\centering
	\includegraphics[width=0.8\linewidth]{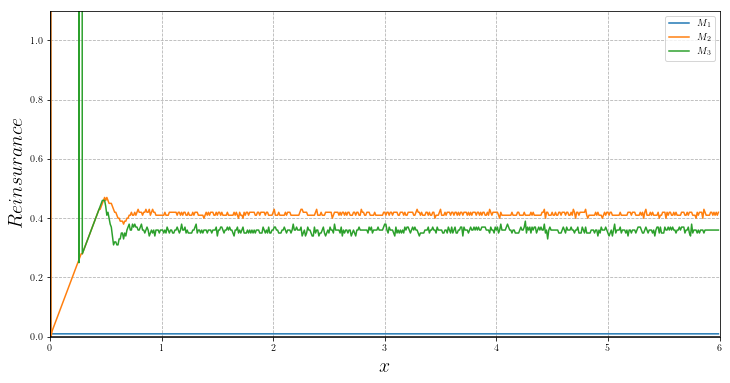}
	\label{XLXLXL}
\end{subfigure}
	\caption{Optimal reinsurance strategies for states $(i)-(viii)$.}
\label{optimalresultone21}
\end{figure}

\begin{figure}
	\centering
	\begin{subfigure}{.4\textwidth}
		\caption{}
		\centering
		\includegraphics[width=.8\linewidth]{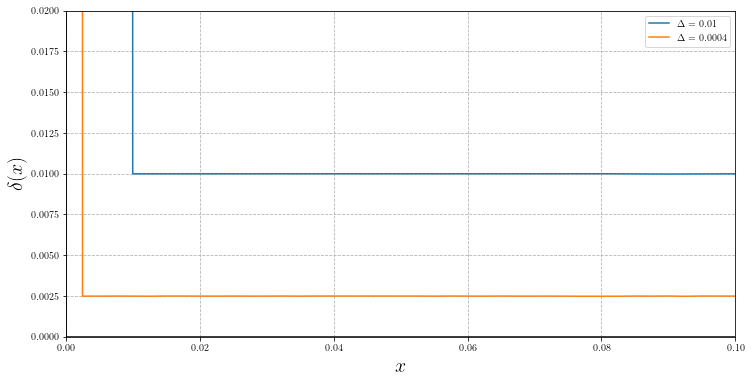}
		\label{XLdelta}
	\end{subfigure}
	\begin{subfigure}{.4\textwidth}
			\caption{}
		\centering
	\includegraphics[width=.8\linewidth]{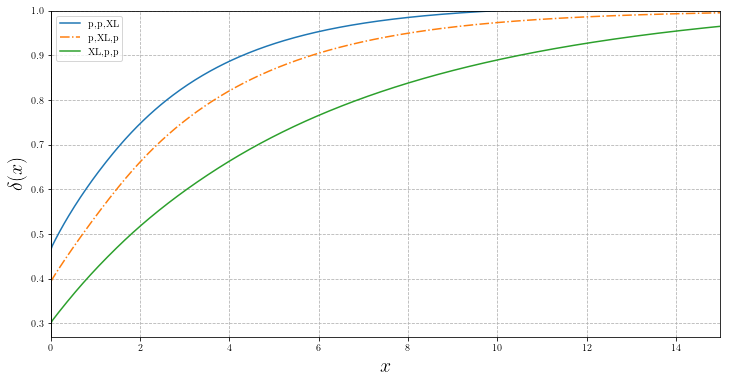}
	\label{oneXL}
	\end{subfigure}\\
	\begin{subfigure}{.4\textwidth}
	\caption{}
	\centering
	\includegraphics[width=.8\linewidth]{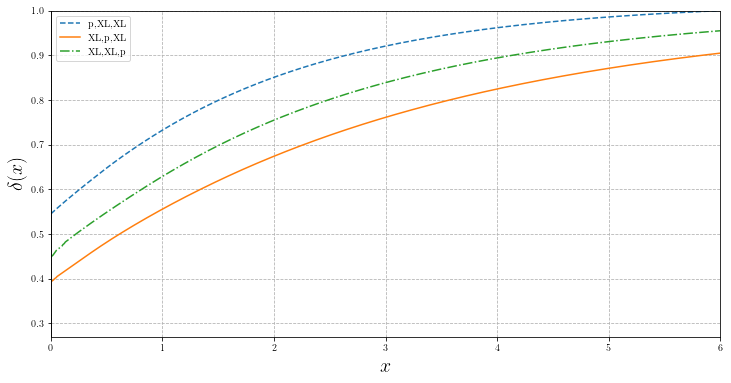}
	\label{twoXL}
\end{subfigure}
	\begin{subfigure}{.4\textwidth}
	\caption{}
		\centering
	\includegraphics[width=.8\linewidth]{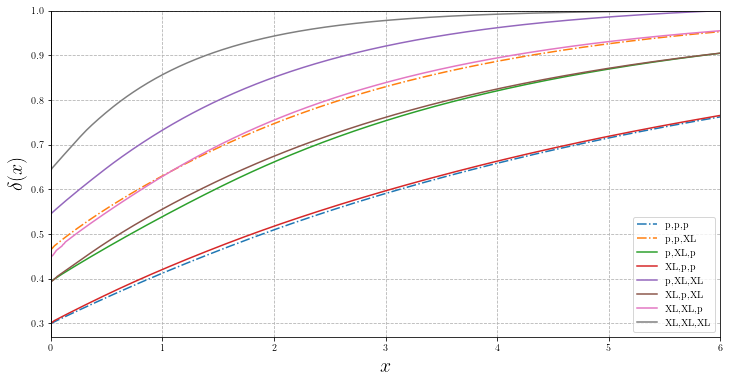}
	\label{allsurvival}
\end{subfigure}
	\caption{(a) The numerical solution of the function $M_1(.)$ for $h=0.01$ and $h=0.0004$ (b) Optimal survival functions for cases where only one XL contract is used (c) Optimal survival functions for cases that use two XL contracts (d) Survival functions for all possible cases}
	\label{optimalresultsurvivalfunctions2}
\end{figure}

Suppose an insurance company is operating on three lines. In the i'th line $i=1,2,3$, the
reinsurance is displayed as $R_i$
and the distribution function is demonstrated as $F_i$. Furthermore, for the claim numbers, Poisson distribution with the parameter $\beta_i$ is used. Suppose further that $R_i$ is one of the proportional reinsurance strategies $(R_p)$ or excess-of-loss $(R_{XL})$. In this condition, there are different states imaginable for choosing the type of reinsurance contract; some of the states will be dealt with later on.
\begin{figure}
	\begin{subfigure}{.4\textwidth}
		\caption{}
		\centering
		\includegraphics[width=1\linewidth]{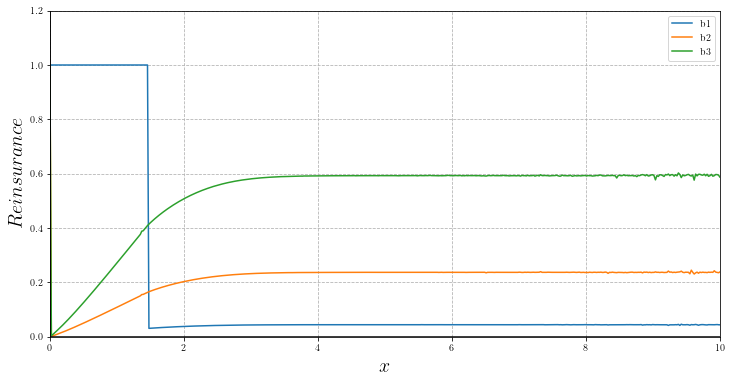}
	\end{subfigure}
\centering
	\begin{subfigure}{.4\textwidth}
				\caption{}
		\includegraphics[width=1\linewidth]{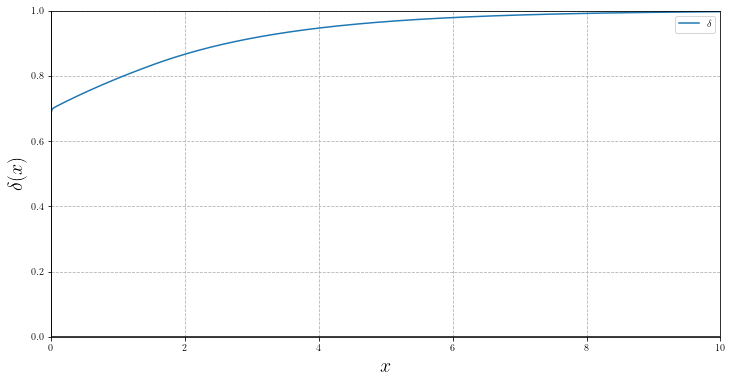}
	\end{subfigure}
	\caption{Three lines of insurance with three proportional reinsurance: (a) Optimal retained Proportion. (b) Survival probability function with Proportional reinsurance. }
	\label{optimalresultproportionmultiple}
\end{figure}
\begin{example}
Let $F_i(x)=1-e^{-\lambda_i x}$, $i=1,2,3$, and $\lambda_1=0.15,\,\,\lambda_2=0.8,\,\, \lambda_3=2, \,\, \beta_1=1, \,\, \beta_2=13,\,\, \beta_3=30, \,\, \eta=3$ and $\eta_1=3.5$. The reinsurance strategy in $i$th line is depicted by $R_i$. As was mentioned before, if $R_i \in \mathcal{R}_p$ then $R_i(y)=b_i(.) y$ and if $R_i \in \mathcal{R}_{XL}$ then $R_i(y)=\min(y,M(.))$, where $b_i(.)$ and $M_i(.)$ are functions of the company's capital.
 If the insurance company considers a reinsurance contract for three lines, the optimization issue will be equal with the uni-dimensional model scrutinized by \cite{azcue2014stochastic}. Figure \ref{optimalresultproportionmultiple}, illustrates rules for the state in which three proportional reinsurance strategies for three lines have been used (this contract could vary from one line to another), followed by  Figure \ref{optimalresultproportiononeone}, displaying the state in which one proportional reinsurance is used for all three lines. Furthermore, If we supplant the proportional reinsurance contract with an excess-of-loss reinsurance, the results would be as Figure \ref{optimalresultMultipleXLRe}, where \ref*{showmultipleXLRe}, illustrates  three reinsurance strategies in three insurance lines, and \ref*{showM2XLRe} and \ref*{showM3XLRe} indicate the charts for reinsurancee strategies of line 2 and 3 (which are the same). In both states, in the sense of increasing the survival function, the use of three appropriate reinsurance contracts would culminate in better results than the use of one contract for all three lines. 
\begin{figure}
	\centering
	\begin{subfigure}{.4\textwidth}
		\caption{}
		\centering
		\includegraphics[width=1\linewidth]{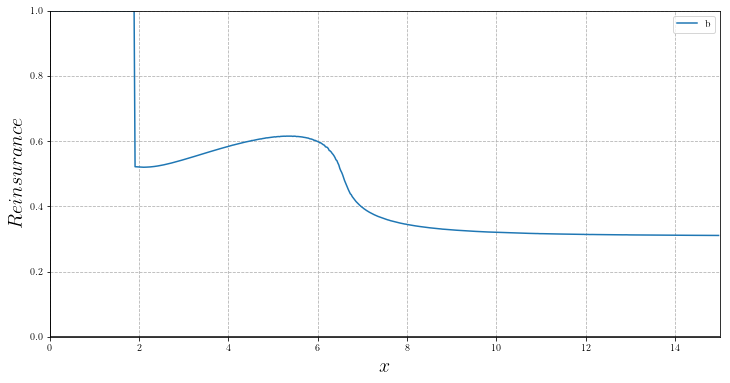}
	\end{subfigure}
	\begin{subfigure}{.4\textwidth}
		\caption{}
		\centering
		\includegraphics[width=1\linewidth]{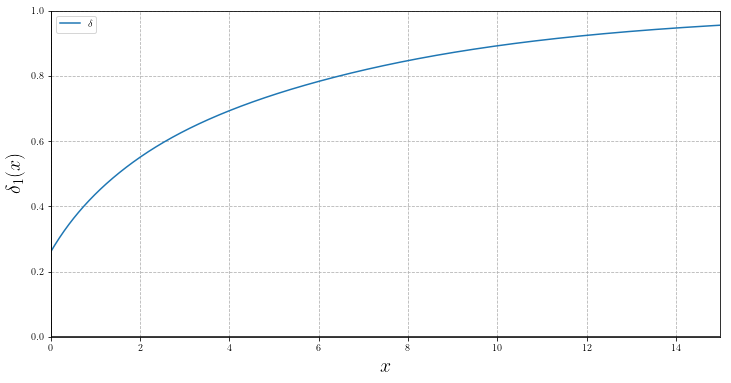}
	\end{subfigure}
	\caption{Three lines of insurance with one proportional reinsurance: (a) Optimal retained Proportion. (b) Survival probability function with reinsurance. }
	\label{optimalresultproportiononeone}
\end{figure}
\begin{figure}
	\centering
	\begin{subfigure}{.4\textwidth}
		\caption{}
		\centering
		\includegraphics[width=1\linewidth]{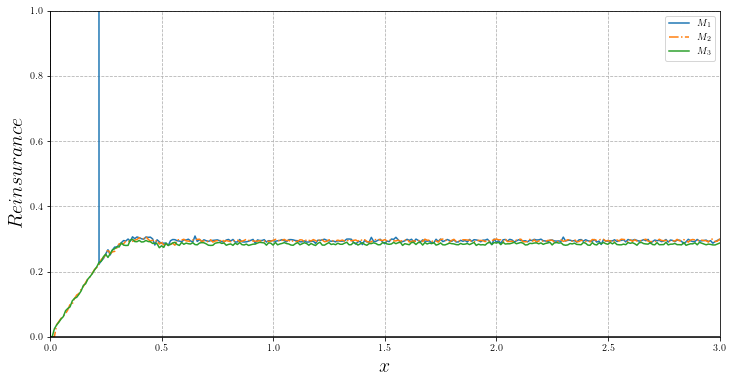}
		\label{showmultipleXLRe}
	\end{subfigure}
	\begin{subfigure}{.4\textwidth}
		\caption{}
		\centering
		\includegraphics[width=1\linewidth]{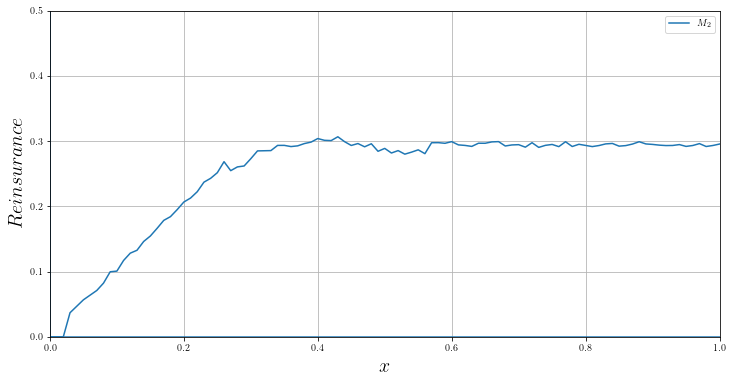}
		\label{showM2XLRe}
	\end{subfigure}
	\\
	\begin{subfigure}{.4\textwidth}
		\caption{}
		\centering
		\includegraphics[width=1\linewidth]{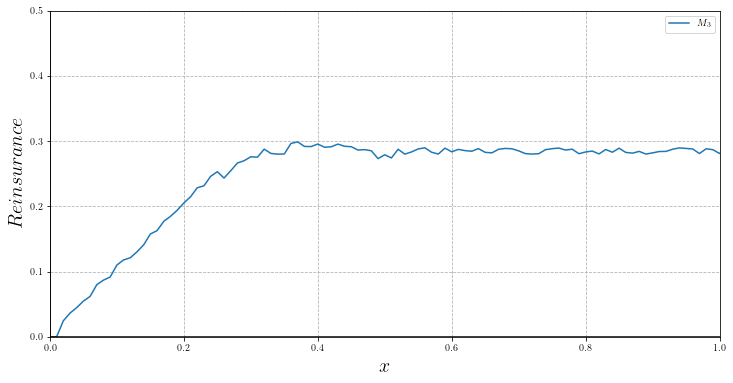}
		\label{showM3XLRe}
	\end{subfigure}
	\begin{subfigure}{.4\textwidth}
		\caption{}
		\centering
		\includegraphics[width=1\linewidth]{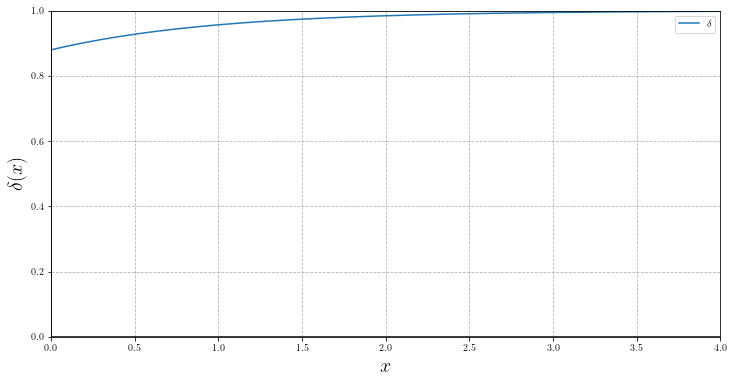}
		\label{SurvivalprobabilityfunctionwithMultipleXLreinsurance}
	\end{subfigure}
	\caption{Three lines of insurance with three excess-of-loss reinsurance: (a) Optimal retention levels for excess-of-loss reinsurances on the first, second and third lines. (b) Optimal retention level for excess-of-loss reinsurance on second line. (c) Optimal retention level for excess-of-loss reinsurance on third line(d) Survival probability function with excess-of-loss  reinsurance.}
	\label{optimalresultMultipleXLRe}
\end{figure}

In Figure \ref{optimalresultsurvivalfunctions}, the diagrams for the optimal  survival function in the states of three excess-of-loss reinsurance strategies for three lines, one  excess-of-loss reinsurance strategy for three lines, three proportional reinsurance strategies for three lines, one proportional reinsurance strategy for three lines  and without reinsurances have been juxtaposed. 
\end{example}
As evident from the figures (\ref{optimalresultproportionmultiple}-\ref{optimalresultsurvivalfunctions}), one could note that these results could be useful for the insurer such that he/she can draw different contracts. By doing this, the insurer can decrease the probability of bankruptcy  \textit{vis-\`{a}-vis } the state that he/she uses only one contract. 
In the previous example, the  exponential distribution is considered for the claims size, and only one type of reinsurance strategy is used to transfer the risk to the secondary insurer on all lines.
In practice, however, the distribution of claims in each line may be completely different from the other one and the insurance company may use different types of reinsurance contracts.
\begin{example}
	In this example, we consider the light-tailed distribution $F_1(x)=1-e^{-0.5 x}$ for the claims size in the first line, the heavy-tailed distribution $F_2(x)=1-(\frac{3}{3+x})^{3}$ for claims size in the second line and the mixture of heavy-tailed distribution and light-tailed distribution $F_3(x)=0.7 F_1(x)+ 0.3 F_2(x)$ for claims size in the last line. 
	If one of the excess-of-loss or proportional reinsurance strategies is used to control the risk of each line, then all possible settings are as follows:
  \begin{multicols}{2}
	\begin{itemize}
		\item[(i)] $R_i\in \mathcal{R}_{p}, \,i=1,2,3$,
	\item[(ii)] $R_i\in \mathcal{R}_{p}, \,i=1,2$ and $R_3\in \mathcal{R}_{XL}$
		\item[(iii)] $R_i\in \mathcal{R}_{p}, \,i=1,3$ and $R_2\in \mathcal{R}_{XL}$
		\item[(iv)] $R_i\in \mathcal{R}_{p}, \,i=2,3$ and $R_1\in \mathcal{R}_{XL}$
	\item[(v)] $R_i\in \mathcal{R}_{XL}, \,i=1,2$ and $R_3\in \mathcal{R}_{p}$
\item[(vi)] $R_i\in \mathcal{R}_{XL}, \,i=1,3$ and $R_2\in \mathcal{R}_{p}$
\item[(vii)] $R_i\in \mathcal{R}_{XL}, \,i=2,3$ and $R_1\in \mathcal{R}_{p}$
	\item[(viii)] $R_i\in \mathcal {R}_{XL}, \,i=1,2,3$.
	\end{itemize}
\end{multicols}
The optimization results for all of the above states are presented in Figure \ref{optimalresultone21}. As it is displayed in the Figure, in all states in which XL has been considered for the line 1, $M_1(x)$ is equal to zero. In Figure \ref*{XLdelta}, for $h=0.01$  and $h=0.0004$,$ M_1(x)$ has been presented for state $(viii)$. In effect, the smaller the $\Delta$  in the FDM, the closer the m1 to zero.  Moreover, in Figure \ref{optimalresultsurvivalfunctions2}, the graph pertained to survival probability is reported in the all 8 above states.  As it is displayed in the Figure, the best state is the one in which all lines resort to  the  XL reinsurance contract. As seen in Figure \ref*{oneXL}, if we decide to resort to the XL strategy only in one line, the best line would be line 2 in which the probability of bigger claims’ size occurring would be higher from the other two lines; furthermore, as it is observed in the Figure \ref{twoXL}, if we opt for using the XL strategy in two lines, lines 2 and 3 would be efficient choices. 
\end{example}
\

The results and findings of this paper reveal that the use of optimal dynamic reinsurance strategies increases the survival function in contrast to the use of one optimal reinsurance strategy in all lines, where a plausible question was raised on the practicality at resorting to the optimal reinsurance. Perhaps at first glance the answer could be a simple no. Indeed, for an insurance company and the reinsurer, implementing a contract which is dependent on the state of the surplus process and at any moment it may change could be extremely difficult, if not impossible. Nonetheless, it can be concluded that there is an interesting point discernible in Figures, indicating the fact that the XL strategy is a more viable option than the proportional strategy in terms of increasing the survival function. Moreover, ostensibly the XL reinsurance contracts oscillate under the circumstances that the capital of the company dwindle and the company is in the critical condition; they seem constant for increased capital. For instance, consider the previous example; if three XL reinsurances for three lines have been used, as long as the company’s assets belong to interval $(0.5 , \infty)$, the reinsurance contracts are constant and they only change when the assets of the company belong to the interval $(0,0.5)$. 
\section{Discussion}
The issue of dynamic reinsurer strategy is an interesting and efficacious approach for augmenting the survival of an insurance company and one can find a great number of studies in literature on the topic of maximizing the survival function with respect to the reinsurance strategy.
 In order to solve this optimization issue, an HJB equation associated with the survival function is adapted. 
 Often, the survival function does not have the smoothness properties needed to interpret it as a solution for the corresponding HJB equation in the classical sense, but it is satisfying in this equation in a weaker concept. The present paper deals with the issue of maximization of the survival function with respect to the dynamic reinurance strategy with this difference that the insurance company shares the potential risk with several reinsurers each pertained to a specified line. The optimal survival function is characterized as the unique non-decreasing viscosity solution of (2.9) with limit one at infinity. 
Unfortunately, obtaining a closed form for the survival function or the reinsurance strategy in the issue discussed at this paper is fiendishly complicated or well-nigh impossible. Therefore, it was more feasible to adopt the numerical solution.  For constructing a numerical solution, the FDM has been employed due to the fact that the convergence of the numerical solution to the survival function can be proved through the techniques prevalent in the literature. The convergent findings are displayed in section 4. 

The results  of the present paper give the insurances companies this stupendous opportunity to share their risk with the reinsurers. In section 5, there are some examples that indicate the fact that using this approach, the survival function will be increased. In a nutshell, with the implementation of this dynamic method for drawing the vector of the reinsurance contracts, the probability of bankruptcy might diminish significantly.

\section*{Appendix}
 Proofs of  Lemma \ref{lemma2.1}, Lemma \ref{S-lipschitz} and Lemma \ref{Sprop2.13} are very similar to Lemma 1.1, Proposition 2.4 and Proposition 2.12 in \cite{azcue2014stochastic},  we therefore omit them.
 
 	{\bf  Proof of Lemma \ref{lemma3.2.'}}. 
 Let us prove first that $H$ is left upper semicontinuous. For assumed $x_0$, $x_k \nearrow x_0$, consider reinsurance strategies ${\boldsymbol{R}}^{(k)}\in \boldsymbol{\mathcal{R}}$ such that
 \begin{align}\label{A.5}
 \sup_{\boldsymbol{R}\in \boldsymbol{\mathcal{R}}} {\tilde{\cal{L}}}_{\boldsymbol{R}}(u)(x_k)\le {{\cal{L}}}_{{\boldsymbol{R}}^{(k)}}(u)(x_k)+\frac{1}{k}.
 \end{align}
 Then, the following is straightforward,
 \begin{align*}
 {{\cal{L}}}_{{\boldsymbol{R}}^{(k)}}(u)(x_0)=&{{\cal{L}}}_{{\boldsymbol{R}}^{(k)}}(u)(x_k)+p_{{\boldsymbol{R}}^{(k)}}(u'(x_0)-u'(x_k))-\beta(u(x_0)-u(x_k))\\
 &+\int_{0}^{x_k}(u(x_0-\alpha)-u(x_k-\alpha)) dF_{{\boldsymbol{R}}}^{(k)}(\alpha)\\
 &+\int_{x_k}^{x_0}u(x_0-\alpha) dF_{{\boldsymbol{R}}}^{(k)}(\alpha)\\
 \ge &{{\cal{L}}}_{{\boldsymbol{R}}^{(k)}}(u)(x_k)+p_{{\boldsymbol{R}}^{(k)}}(u'(x_0)-u'(x_k))-\beta(u(x_0)-u(x_k))\\
 &+\int_{0}^{x_k}(u(x_0-\alpha)-u(x_k-\alpha)) dF_{{{\boldsymbol{R}}}^{(k)}}(\alpha).
 \end{align*} 
 So, we have,
 \begin{align}\label{A.6}
 \limsup_{k\rightarrow \infty}{{\cal{L}}}_{{\boldsymbol{R}}^{(k)}}(u)(x_0)\ge \limsup_{k\rightarrow \infty}{{\cal{L}}}_{{\boldsymbol{R}}^{(k)}}(u)(x_k).
 \end{align}
 Then, from (\ref{A.5}) and (\ref{A.6}), the following result can be derived; 
 $$\sup_{\boldsymbol{R}\in \boldsymbol{\mathcal{R}}} {{\cal{L}}}_{\boldsymbol{R}}(u)(x_0)\ge \limsup_{k\rightarrow \infty}{{\cal{L}}}_{{\boldsymbol{R}}^{(k)}}(u)(x_0)\ge \limsup_{k\rightarrow \infty}\left(\sup_{\boldsymbol{R}\in \boldsymbol{\mathcal{R}}} {{\cal{L}}}_{\boldsymbol{R}}(u)(x_k)\right).$$
 Consequently, the following relation is dominant
 $\limsup_{x\rightarrow x_0^-}\left(\sup_{\boldsymbol{R}\in \boldsymbol{\mathcal{R}}} {{\cal{L}}}_{\boldsymbol{R}}(u)(x)\right)\le \sup_{\boldsymbol{R}\in \boldsymbol{\mathcal{R}}} {{\cal{L}}}_{\boldsymbol{R}}(u)(x_0).$
 Now, we must prove the following    $\limsup_{x\rightarrow x_0^+}\left(\sup_{\boldsymbol{R}\in \boldsymbol{\mathcal{R}}} {{\cal{L}}}_{\boldsymbol{R}}(u)(x)\right)\le \sup_{\boldsymbol{R}\in \boldsymbol{\mathcal{R}}} {{\cal{L}}}_{\boldsymbol{R}}(u)(x_0).$
 Given any sequence $x_k\searrow x_0$, take reinsurance strategies ${\boldsymbol{R}}^{(k)}\in \boldsymbol{\mathcal{R}}$ such that 
$$
 \sup_{\boldsymbol{R}\in \boldsymbol{\mathcal{R}}} {{\cal{L}}}_{\boldsymbol{R}}(u)(x_k)\le {{\cal{L}}}_{{\boldsymbol{R}}^{(k)}}(u)(x_k)+\frac{1}{k}.
$$
 If one of the reinsurance contracts is LXL reinsurance, for example $\mathcal{R}_1=\mathcal{R}_{LXL}$, 
 take $\bar{{\boldsymbol{R}}}^{(k)}=(\bar{R}_1^{(k)},\bar{R}_2^{(k)},\cdots , \bar{R}_n^{(k)}) \in \boldsymbol{\mathcal{R}}$ such that
 \begin{align*}
 \bar{R}_1^{(k)}(\alpha)=
 \begin{cases}	R_1^{(k)}(\alpha)	& if \,\, R_1^{(k)}(\alpha)=\alpha\,\, for\,\, all\,\, \alpha \\
 R_1^{(k)}(\alpha)	& if \,\, R_1^{(k)}(\alpha)=a_k\wedge \alpha+(\alpha-L-a_k)^+\,\, with\,\, a_k \notin (x_0,x_k) \\
 \alpha \wedge x_0 +(\alpha-L-x_0)^+	& if \,\, R_1^{(k)}(\alpha)=a_k\wedge \alpha+(\alpha-L-a_k)^+\,\, with\,\, a_k \in (x_0,x_k) 
 \end{cases}
 \end{align*}
 and $\bar{R}_2^{(k)}=R_2^{(k)}, \cdots, \bar{R}_n^{(k)}=R_n^{(k)}$.
 If
 $a_k \le x_0$ 
 then
 \begin{align*}
 &\int_{0}^{\infty} (u(x_k-\alpha)-u(x_0-\alpha))dF_{\bar{{\boldsymbol{R}}}^{(k)}}(\alpha)\\
 &= \int_{0}^{\infty} (u(x_k-\alpha)-u(x_0-\alpha))dF_{{\boldsymbol{R}}^{(k)}}(\alpha)\\
 &\le \sup_{x \in [0,A]}|u'(x)|(y-x_0).
 \end{align*}
 Let us define $ {{\boldsymbol{R}}}_{n-1}^{(k)}=(R_2^{(k)},\cdots, R_n^{(k)})$. If
 $a_x \in (x_0,x_k)$, 
 then
 \begin{align*}
 &\int_{0}^{\infty}  u(x_k-\alpha)dF_{{{\boldsymbol{R}}}^{(k)}}(\alpha)- \int_{0}^{\infty} u(x_0-\alpha)dF_{\bar{{\boldsymbol{R}}}^{(k)}}(\alpha)\\
 =&
 \int_{0}^{\infty} \int_{0}^{\infty} u(x_k-R_1^{(k)}(\alpha_1)-\alpha_2)dF_{{{\boldsymbol{R}}}_{n-1}^{(k)}}(\alpha_2)dF_1(\alpha_1)-\int_{0}^{\infty} \int_{0}^{\infty} u(x_0-\bar{R}_1^{(k)}(\alpha_1)-\alpha_2)dF_{\bar{{\boldsymbol{R}}}_{n-1}^{(k)}}(\alpha_2)dF_1(\alpha_1)\\
 =&\int_{0}^{\infty} \int_{0}^{\infty} u(x_k-R_1^{(k)}(\alpha_1)-\alpha_2)dF_{{{\boldsymbol{R}}}_{n-1}^{(k)}}(\alpha_2)dF_1(\alpha_1)-\int_{0}^{\infty} \int_{0}^{\infty} u(x_0-\bar{R}_1^{(k)}(\alpha_1)-\alpha_2)dF_{{{\boldsymbol{R}}}_{n-1}^{(k)}}(\alpha_2)dF_1(\alpha_1)\\
 =& \int_{0}^{\infty} \int_{0}^{\infty}\big( u(x_k-R_1^{(k)}(\alpha_1)-\alpha_2)-u(x_0-\bar{R}_1^{(k)}(\alpha_1)-\alpha_2)dF_{{{\boldsymbol{R}}}_{n-1}^{(k)}}(\alpha_2)\big)dF_1(\alpha_1)\\
 =& \int_{0}^{x_0} \int_{0}^{\infty}\big( u(x_k-R_1^{(k)}(\alpha_1)-\alpha_2)-u(x_0-R_1^{(k)}(\alpha_1)-\alpha_2)dF_{{{\boldsymbol{R}}}_{n-1}^{(k)}}(\alpha_2)\big)dF_1(\alpha_1)\\
 &+ \int_{x_0}^{a_k} \int_{0}^{\infty}\big( u(x_k-\alpha_1-(\alpha_1-L-a_k)^+-\alpha_2)-u(-(\alpha_1-L-x_0)^+-\alpha_2)dF_{{{\boldsymbol{R}}}_{n-1}^{(k)}}(\alpha_2)\big)dF_1(\alpha_1)\\
 &+ \int_{a_k}^{x_k} \int_{0}^{\infty}\big( u(x_k-\alpha_k-(\alpha_1-L-a_k)^+-\alpha_2)-u(-(\alpha_1-L-x_0)^+-\alpha_2)dF_{{{\boldsymbol{R}}}_{n-1}^{(k)}}(\alpha_2)\big)dF_1(\alpha_1)\\
 &+ \int_{x_k}^{\infty} \int_{0}^{\infty}\big( u(x_k-\alpha_k-(\alpha_1-L-a_k)^+-\alpha_2)-u(-(\alpha_1-L-x_0)^+-\alpha_2)dF_{{{\boldsymbol{R}}}_{n-1}^{(k)}}(\alpha_2)\big)dF_1(\alpha_1)\\
 &\le \sup_{x \in [0,A]}|u'(x)|(x_k-x_0)\\
 &+\sup_{\alpha \in [x_0,a_k]}(u(x_k-\alpha-(\alpha-L-a_k)^+)-u(-(\alpha-L-x_0)^+))p(x_0\le \alpha \le a_k)\\
 &+\sup_{\alpha \in [a_k,x_k]}(u(x_k-a_k-(\alpha-L-a_k)^+)-u(-(\alpha-L-x_0)^+))p(a_k\le \alpha \le x_k)\\
 &+\sup_{\alpha \in [x_k,\infty)}(u(x_k-a_k-(\alpha-L-a_k)^+)-u(-(\alpha-L-x_0)^+))p( \alpha \ge x_k)
 \end{align*}
 According to the property of  right-continuously of  distribution function; if 
 $x_k\searrow x_0$,
 then
$$p(x_0\le \alpha \le a_k)\longrightarrow 0\,\,\,\,
 \text{and}\,\,\,\,
 p(a_k\le \alpha \le x_k)\longrightarrow 0.$$
 Now, the term 
 $\sup_{\alpha \in [x_k,\infty)}(u(x_k-a_k-(\alpha-L-a_k)^+)-u(-(\alpha-L-x_0)^+))$ 
 should become the focus of attention. In this case, there are two situations as outlined below:
 \begin{itemize}
 	\item [(I)]If, there is a finite value $ m$ satisfying the following,
 	\begin{align*}
 	&\sup_{\alpha \in [x_k,\infty)}(u(x_k-a_k-(\alpha-L-a_k)^+)-u(-(\alpha-L-x_0)^+))\\
 	&=u(x_k-a_k-(m-L-a_k)^+)-u(-(m-L-x_0)^+).
 	\end{align*}
 	So, 
 	\begin{align*}
 	& \lim_{x_k\searrow x_0}\sup_{\alpha \in [x_k,\infty)}(u(x_k-a_k-(\alpha-L-a_k)^+)-u(-(\alpha-L-x_0)^+))\\
 	&=\lim_{x_k\searrow x_0}(u(x_k-a_k-(m-L-a_k)^+)-u(-(m-L-x_0)^+))\\
 	&=u(-(m-L-x_0)^+)-u(-(m-L-x_0)^+)=0.
 	\end{align*}  
 	\item[(II)]  
 	If
 	\begin{align*}
 	& \sup_{\alpha \in [x_k,\infty)}(u(xk-a_k-(\alpha-L-a_k)^+)-u(-(\alpha-L-x_0)^+))\\
 	&=\lim_{\alpha \rightarrow \infty}u(x_k-a_k-(\alpha-L-a_k)^+)-u(-(\alpha-L-x_0)^+)
 	\end{align*}
 	then
 $$
 	\sup_{\alpha \in [x_k,\infty)}(u(xk-a_k-(\alpha-L-a_k)^+)-u(-(\alpha-L-x_0)^+))
 	=u(x_k-a_k)-u(0).
 $$
 	So,
 	\begin{align*}
 	&\lim_{x_k\searrow x_0}\sup_{\alpha \in [x_k,\infty)}(u(x_k-a_k-(\alpha-L-a_k)^+)-u(-(\alpha-L-x_0)^+))\\
 	&=\lim_{x_k\searrow x_0}(u(x_k-a_k)-u(0))=0.
 	\end{align*}  
 \end{itemize}
 It should be noted that  
 $R_1^{(k)}(\alpha)-(x_k-x_0)\le \bar{R}_1^{(k)}(\alpha)\le R_1^{(k)}(\alpha).$
 So 
 $p_{\bar{R}_1^{(k)}} \nearrow p_{R_1^{(k)}}$. 
 Thus
 \begin{align*}
 H(x_k)-H(x_0)&\le {{\cal{L}}}_{\boldsymbol{R}^{(k)}}(u)(x_k)-{{\cal{L}}}_{{\boldsymbol{R}}^{(k)}}(u)(x_0)+\epsilon\\
 &\le p_{\boldsymbol{R}^{(k)}}\sup_{x\in [0,A]}|u''(x)|(x_k-x_0)+2\beta \sup_{x\in [0,A]}|u'(x)|(x_k-x_0)\\
 &+\sup_{\alpha \in [x_0,a_k]}(u(x_k-\alpha-(\alpha-L-a_k)^+)-u(-(\alpha-L-x_0)^+))p(x_0\le \alpha \le a_k)\\
 &+\sup_{\alpha \in [a_k,x_k]}(u(x_k-a_k-(\alpha-L-a_k)^+)-u(-(\alpha-L-x_0)^+))p(a_k\le \alpha \le x_k)\\
 &+\sup_{\alpha \in [x_k,\infty)}(u(x_k-a_k-(\alpha-L-a_k)^+)-u(-(\alpha-L-x_0)^+))p( \alpha \ge x_k)
 \end{align*}
 and so we get that $H$ is right upper semicontinuous. 
 The proof for the case $\mathcal{R}_1=\mathcal{R}_{XL}$ and $\mathcal{R}_1=\mathcal{R}_{p}$ are simpler,  we therefore omit them.
 Now, repeating the arguments presented in the proof of the Lemma 3.2 of \cite{azcue2014stochastic} (replacing   $\boldsymbol{R}=(\boldsymbol{R}_t)_{t \ge 0}=({R_1}_t,\cdots, {R_n}_t)_{t \ge 0}$ with  $\bar{R}$), the following relation is obtained,
 $$|H(y)-H(x_0)|\le K_{A}(y-x_0+F(y)-F(x_0)).$$
 The proof is complete.
 $\square$\\
 
If in  Lemma 4.2 and Lemma 4.3 of \cite{azcue2014stochastic}, replace 
$${\cal{L}}_{R}(\delta)(x)=p_{R}\delta'(x)-\beta\delta(x)+\beta\int_{0}^{x}\delta(x-\alpha)dF_{R}(\alpha)$$
 with 
 $${\cal{L}}_{\boldsymbol{R}}(\delta)(x)=p_{\boldsymbol{R}}\delta'(x)-(\sum_{i=1}^{n}\beta_i)\delta(x)+(\sum_{i=1}^{n}\beta_i)\int_{0}^{x}\delta(x-\alpha)dF_{\boldsymbol{R}}(\alpha),$$ then  two lemmas Lemma \ref{lemm4.2} and \ref{lemm4.3} are derived and
to prove the Lemma \ref{Lemma4.4}, it suffices to repeat the proof of Lemma 4.4 of \cite{azcue2014stochastic} by replacing   $\boldsymbol{R}=(\boldsymbol{R}_t)_{t \ge 0}=({R_1}_t,\cdots, {R_n}_t)_{t \ge 0}$ with  $\bar{R}$. Therefore,  proofs of these Lemmas are omitted. 

{\bf proof of Lemma \ref{numeric-increasing}} 
(i)  Assume that $i$ is a positive integer with $f '_h(k h)\ge  0,  k = 1, ..., i- 1$. Then 
$f_h(0) \le f_h(h) \le \cdots \le f_h((i-1)h)$ and thus 
$$G_{\boldsymbol{R}}(ih)=\sum_{\{j\le i\}}f_h((i-j)h)P\{(j-1)h<\boldsymbol{R}(Y)\le jh\}\le f_h((i-1) h).$$
 So for  $s=ih$
$$\frac{\beta f_h(s-h)-\beta G_{\boldsymbol{R}}(s)}{p_{{\boldsymbol{R}}}} \ge 0,$$
and thus obviously that
$$f'_h(ih)=f'_h(s)=\inf_{\boldsymbol{R}\in \mathcal{R}}\frac{\beta f_h(s-h)-\beta G_{\boldsymbol{R}}(s)}{p_{{\boldsymbol{R}}}} \ge 0,$$
which completes the induction.

(ii) The proof of this part is similar to the proof of Lemma 9 in \cite{nozadi2014optimal}, section 3.3.
$\square$

\bibliographystyle{agsm} 
\bibliography{references.bib}
\end{document}